\documentclass[draft]{article}
\usepackage[english]{babel}
\usepackage{amssymb}
\usepackage{amsfonts}
\usepackage{amsthm}

\renewcommand{\emptyset}{{\mathchar28}}

\newcommand{\bbC}{\mathord{\mathbb C}}
\newcommand{\bbR}{\mathord{\mathbb R}}
\newcommand{\Tr}{\mathop{\rm Tr}\nolimits}
\newcommand{\ad}{\mathop{\rm ad}\nolimits}
\newcommand{\Ad}{\mathop{\rm Ad}\nolimits}

\newcommand{\e}{\mathrm e}
\newcommand{\const}{{\mathop{\rm const}}}
\newcommand{\one}{{\mathbf1}}

\newtheorem{lemma}{Lemma}
\newtheorem{proposition}{Proposition}
\newtheorem{theorem}{Theorem}
\newtheorem{corollary}{Corollary}
\newtheorem{remark}{Remark}

\title{On complete affine structures in Lie groups \footnote{This
is a translation of the paper published in book: Essays in
differential equations and algebra, IITPM, Omsk (1992), 67-80 (in
Russian). There are nonessential changes in the exposition,
references are extended but the mathematical content is identical.
In the case of the field $\bbC$, Section~1 overlaps with the paper
\cite{Hel}.}}
\author{V.M.\,Gichev}
\date{}
\begin{document}
\maketitle
\begin{abstract}
Left invariant affine structures in a Lie group $G$ are in
one-to-one correspondence with left-symmetric algebras over its
Lie algebra $\mathfrak g=T_eG$ (``over'' means that the commutator
$[x,y]=xy-yx$ coincides with the Lie brackets; left-symmetric
algebras can be defined as Lie-admissible algebras such that the
multiplication by left defines a representation of the underlying
Lie algebra). An affine structure (and the corresponding left
symmetric algebra) is complete if $G$ is affinely equivalent to
$\mathfrak g$. By the main result of this paper, a complete left
symmetric algebra admits a canonical decomposition: there is a
Cartan subalgebra $\mathfrak h$ such that the root subspaces for
the representations $L$ (by left multiplications) and $\ad$
coincide. Then operators $L(x)$ and $\ad(x)$ have equal semisimple
parts for all $x\in\mathfrak h$. This decomposition is unique. For
simple complete left-symmetric algebras whose canonical
decomposition consists of one dimensional spaces we define two
types of graphs and prove some their properties.  This makes
possible to describe, for dimensions less or equal to $5$, these
graphs and algebras.
\end{abstract}

A connection in a manifold $G$ can be defined by the covariant
derivative $\nabla_\xi\eta$, where $\xi,\eta$, and
$\nabla_\xi\eta$, are smooth vector fields. If the torsion and the
curvature vanish
\begin{eqnarray*}
&\nabla_{\xi} \eta -\nabla_{\eta} \xi
-[\xi,\eta] =0,\\
&\nabla_{\xi}\nabla_{\eta}\zeta - \nabla_{\eta}\nabla_{\xi} \zeta
-\nabla_{[\xi,\eta]} \zeta =0,
\end{eqnarray*}
where $[\xi,\eta]$ is the Lie bracket of vector fields, then the
manifold admits an atlas of coordinate charts with affine
coordinate transformations. This defines the {\it affine
structure} in $G$. If $G$ is a Lie group and $\nabla, \xi, \eta$
are left invariant then the field  $\nabla_{\xi} \eta$ is also
left invariant. Realizing the Lie algebra $\frak g$ by these
vector fields we get an operation in $\frak g$ which satisfies
identities
\begin{eqnarray}
& \xi\eta - \eta\xi = [\xi,\eta], \label{torsi}\\
&[\xi,\eta] \zeta = \xi(\eta\zeta) - \eta (\xi\zeta).\label{curvi}
\end{eqnarray}
The second can be written as follows:
\begin{eqnarray}
 (\xi,\eta,\zeta)=(\eta,\xi,\zeta),\label{lefts}
\end{eqnarray}
where $(\xi,\eta,\zeta) = \xi(\eta\zeta) - (\xi\eta)\zeta$ is the
{\it associator}. Algebras satisfying (\ref{lefts}) are called
{\it left-symmetric}. They were defined by Vinberg\footnote{and by
many other authors in various forms, see a recent survey
\cite{Bu05}.} in the paper \cite{Vin} which contains the
classification of homogeneous convex pointed cones. Left-symmetric
algebras are {\it Lie--admissible}: (\ref{lefts}) implies that the
brackets in (\ref{torsi}) satisfies the Jacobi identity. By
(\ref{curvi}),
\begin{eqnarray}\label{reprl}
L([\xi,\eta]) = [L(\xi),L(\eta)],
\end{eqnarray}
where $L(\xi)\zeta = \xi\zeta$. Thus left symmetric algebras could
be defined as Lie--admissible algebras such that the left
multiplication is a representation of the underlying Lie algebra.
There is a natural realization of a left-symmetric algebra $\frak
g$ by affine vector fields on the linear space $\frak g$: each $x
\in \mathfrak g$ corresponds the field
\begin{eqnarray}\label{famfi}
F_x(y) =xy+x. 
\end{eqnarray}
For the multiplication defined by the linear part of
(\ref{famfi}), identities (\ref{torsi}) and (\ref{curvi})  are
equivalent to the condition that the set
\begin{eqnarray}\label{afffi}
 \{F_x\colon \; x\in \mathfrak g \}
\end{eqnarray}
is a Lie algebra with respect to the Lie bracket of vector fields.
Let $G$ be the simply connected Lie group  corresponding to
$\mathfrak g$. The realization above defines an action of $G$ by
affine transformations of $\mathfrak g$. Since the mapping $x\to
F_x(o)$ is nondegenerate, the stable subgroup of the origin $o$ is
discrete. This defines a natural affine structure in $G$. If the
action is transitive on $\mathfrak g$ then it is free since the
stable subgroup is trivial being isomorphic to the fundamental
group of the orbit. Then $G$ is diffeomorphic to Euclidean space.
Therefore, $G$ is solvable in this case (note that $G$ cannot be
isomorphic to the universal covering group of $\mathop{\mathrm
SL}(2,\bbR)$ being linear). We say that the affine structure and
the left-symmetric algebra are {\it complete} if the action above
is transitive. This is equivalent to the geodesical completeness.

The paper is organized as follows. In Section~1, we give the
algebraic criterion of the completeness for real left-symmetric
algebras (for complex ones it is known, see \cite{Hel}). Section~2
contains the main result of the paper. We prove that a complete
left-symmetric algebra admits a {\it canonical decomposition}:
there exists a Cartan subalgebra $\mathfrak h$ of $\mathfrak g$
such that the root decompositions for the representations $L$ and
$\ad$ coincide as well as semisimple parts of operators $L(x)$ and
$\ad(x)$ for all $x\in\mathfrak h$. The decomposition is unique.
In the third section, we consider simple left-symmetric algebras
that admit one dimensional root decompositions. We define two
types of graphs for them and prove some their properties. This
makes possible to describe the graphs and give the complete list
of these algebras in dimensions less than $6$~\footnote{a
computation shows that in dimension 6 there is a lot of graphs
with these properties as well as simple algebras of this type.}.

All algebras are finite dimensional, real or complex. We denote
the field by  $\mathbb F$.  In most cases, the field is not
specified since arguments hold for both cases. By $\mathfrak g$ we
denote the left-symmetric algebra as well as the underlying Lie
algebra, by $G$ the corresponding Lie group of affine
transformations of $\mathfrak g$.

There is one more equivalent version  of (\ref{lefts}):
\begin{eqnarray}\label{lrcpm}
[L(\xi),R(\zeta)] = R(\xi\zeta) - R(\zeta)R(\xi),
\end{eqnarray}
where $R(\xi)\eta = \eta\xi$ is the operator of the right
multiplication by $\xi$.

\section{Algebraic criterion of completeness}

Any left-symmetric algebra $\mathfrak g$ can be extended by the
addition of the identity element $\one$ which satisfies $\one\cdot
x = x \cdot \one = x$ for all $x\in \mathfrak g$. A simple
straightforward calculation shows that the extended algebra
$\mathfrak g_{\one } = \mathbb F\cdot \one + \mathfrak g$
satisfies (\ref{lefts}). In geometric terms, this means that the
affine action in the vector space $\mathfrak g$ naturally defines
a linear representation of $G$ in the vector space $\mathfrak
g_{\one }$: any affine transformation $x\to Ax + b$ of $\mathfrak
g$ in the hyperplane $\one + \mathfrak g$ of vectors $\one+x$,
$x\in\mathfrak g$,  uniquely  extends to a linear transformation
of $\mathfrak g_{\one}$. Vector fields (\ref{famfi}) on $\one +
\mathfrak g$ have the form $ F_x(y) = x(\one + y) $ and are
included to the family
\begin{eqnarray}\label{famve}
 \{ F_x\colon\; x\in \mathfrak
g_{\one }\}, 
\end{eqnarray}
where $ F_x(y) = xy$, $y\in \mathfrak g_{\one }$. On the level of
Lie algebras, $\mathfrak g_{\one}$ is the trivial one dimensional
extension of $\mathfrak g$. Let us equip with the index ${\one }$
spaces and groups relating to $\mathfrak g_{\one }$; for instance,
$G_{\one}$ is the linear group in the linear space $\mathfrak
g_{\one }$ corresponding to (\ref{famve}). Set
\begin{eqnarray}
 P_{\one}(x) =\det(R_{\one} (x)), \nonumber \\
 P(x) = \det(I + R(x)), \label{polip} 
\end{eqnarray}
where $I$ is the identical transformation. The polynomial $P$ is
the restriction of $P_\one$ to the hyperplane $\one +\mathfrak g$.
\begin{lemma}\label{eigpe}
For all  $x,y\in \mathfrak g_{\one }$
\begin{eqnarray}
R(\e^{L(y)}x) = \e^{L(y)}R(x)\e^{-\ad (y)}. \label{ripra} 
\end{eqnarray}
\end{lemma}
\begin{proof}
Differentiating  (\ref{ripra}) on $y$ at $y=0$ we get
\begin{eqnarray}\label{difft}
R(L(u)x) = L(u)R(x)-R(x)\ad(u),\quad u\in \mathfrak g_{\one }.
\end{eqnarray}
This is equivalent to the equality
$$
z(ux) = u(zx) - [u,z]x
$$
for all  $u,x,z\in \mathfrak g_{\one }$ that is the same as
(\ref{reprl}). Since left and right sides of  (\ref{ripra}) define
representations of $G_{\one }$ and the equality holds for $y=0$,
(\ref{ripra}) follows from (\ref{difft}).
\end{proof}
\begin{corollary}
For all  $x,y\in \mathfrak g_{\one }$
\begin{eqnarray}
 P_{\one} (\e^{L(y)}x) = \e^{\Tr(R(y))} P_{\one} (x), \label{eigep} 
\end{eqnarray}
\end{corollary}
\begin{proof}
Since $\ad(y)=L(y)-R(y)$, the right side of (\ref{eigep}) is the
determinant of the right side in (\ref{ripra}).
\end{proof}
\noindent According to (\ref{eigep}), the polynomial $P_{\one}$ is
an eigenfunction of $G_{ \one}$:
\begin{eqnarray}
 P_{ \one} (g(x)) = \lambda(g) P_{\one} (x), \label{eigef} 
\end{eqnarray}
where $\lambda$ is some one dimensional character of $G_{ \one }$.
For $x= \one$, $g=\exp(y)$, where $\exp:\frak g_\one\to G_\one$ is
the exponential mapping, setting  $\e^y = \exp(y)( \one)$ we get
\begin{eqnarray*}
 \lambda (\exp(y)) =
P_{ \one} (\e^y) = \e^{\Tr(R(y))} 
\end{eqnarray*}
(note that $P_{ \one}( \one) = 1$).
\begin{lemma}\label{stand}
If $\mathfrak g$ is complete then the operator $L(y)$ is
degenerate for all  $y\in\mathfrak g$.
\end{lemma}
\begin{proof}
Otherwise, the operator $\e^{L(y)}-I$ is nondegenerate for generic
$y\in\mathfrak g$.  Then any affine mapping with the linear part
$\e^{L(y)}$ has a fixed point contradictory to the assumption that
$\frak g$ is complete since transitive actions are free.
\end{proof}
\begin{proposition}\label{realc}
The left-symmetric algebra $\mathfrak g$ is complete if and only
if any of the following conditions holds for all $x\in\mathfrak
g$:
\begin{itemize}
\item[\rm (a)] $R(x)$ is nilpotent;
 \item[\rm (b)] $P(x)=1$;
\item[\rm (c)] $P(x) \ne 0$;
 \item[\rm (d)] $\Tr R(x) =0$.
\end{itemize}
\end{proposition}
\begin{proof}
Implications  (a) $\Rightarrow$ (b) $\Rightarrow$  (c), \quad (a)
$\Rightarrow$ (d) are evident. It follows from (\ref{eigep}) that
(b) and (d) are equivalent. If $\det (I+tR(x))=1$ for all
$t\in\mathbb F$ then $R(x)$ is nilpotent; thus  (a) follows from
(b). Hence it is sufficient to prove that (c) implies the
completeness of $\mathfrak g$ and that (d) is true if $\frak g$ is
complete.

Since vectors of vector fields (\ref{afffi})  at $y$ (\ref{famfi})
have the form $xy + x = (I+R(y))x$, the condition  (c) is
equivalent to the assumption that the tangent space to the orbit
of $y$ is equal to $\mathfrak g$. Therefore, all orbits are open.
Then they are closed; hence there is only one orbit.

Let $\mathfrak g$ be complete, $y\in\frak g$, $\mathfrak g^0$ be
the root subspace of the eigenvalue 0 of $L(y)$ and $\mathfrak
g^1$ be the sum of all other root subspaces. Then $\mathfrak
g_{\one } = \mathfrak g^0 \oplus \mathfrak g^1$. Since  $\mathfrak
g\cdot\mathfrak g_{\one} \subseteq\mathfrak g$, the inclusion
$\mathfrak g^1 \subseteq \mathfrak g$ holds. Therefore,
$$
\mathfrak g^0 \cap(\one + \mathfrak g) \neq\emptyset.
$$
Hence there exist $x\in \one + \mathfrak g$ and $n\in\mathbb N$
such that $L^n(y)x =0$. By (\ref{eigep}),
\begin{eqnarray}\label{pexpl}
P_{\one} (\e^{tL(y)}x) = \e^{t\Tr R(y)} P_{\one}(x) 
\end{eqnarray}
for all $t\in\mathbb F$. The left side of (\ref{pexpl}) is a
polynomial but the right one is an exponent on $t$. Hence either
$P_{\one}(x)=0$ or $\e^{t\Tr R(y)} = \const$. The assumption
$P_{\one}(x) = 0$ contradicts to the completeness of $\mathfrak g$
and (\ref{eigef}). This proves the proposition.
\end{proof}
\begin{remark}\rm
If $\mathbb F=\mathbb C$ and $\frak g$ is complete then (d) is an
evident consequence of (c) since any nonconstant polynomial has
zeroes.
\end{remark}
\noindent Corollaries below follow from Proposition~\ref{realc},
(a) and (d), respectively.
\begin{corollary}\label{compl}
Subalgebras and quotient algebras of a complete left-symmetric
algebra are complete.\qed
\end{corollary}
\begin{corollary}
A real left-symmetric algebra $\frak g$ is complete if and only if
$\frak g\otimes\mathbb C$ is complete.\qed
\end{corollary}

\section{The canonical decomposition}

Let $\mathfrak h$ be a Cartan subalgebra of the Lie algebra
$\mathfrak g$. If $\mathfrak g$ is complex and $\lambda$ is a
linear functional on $\mathfrak h$ then $\mathfrak h^{\lambda}$,
$\mathfrak g^{\lambda}$ denote root subspaces of $\lambda$ for
representations $\ad$, $L$, respectively. For real $\mathfrak g$,
$\lambda$ is a linear functional on $\frak h^\mathbb C=\frak
h+i\frak h$ and $\mathfrak h^{\lambda}$, $\mathfrak g^{\lambda}$
are real parts of the corresponding root spaces.  Then
\begin{eqnarray}
\mathfrak g = \sum_{\lambda}\oplus\, \mathfrak h^{\lambda},\label{decad}\\
\mathfrak g = \sum_{\lambda}\oplus\, \mathfrak
g^{\lambda}\label{decle}.
\end{eqnarray}
We say that the decomposition is {\it canonical} if (\ref{decad})
coincides with (\ref{decle}) (i.\,e.\ $\mathfrak h^{\lambda} =
\mathfrak g^{\lambda}$ for all $\lambda$). The space $\mathfrak
h_\one=\mathbb F\cdot\one+\frak h$ is the Cartan subalgebra for
$\mathfrak g_{\one}$ (and the trivial one dimensional extension of
$\frak h$). Similar decompositions hold for $\mathfrak h_\one$ and
$\mathfrak g_{\one }$:
\begin{eqnarray}
\mathfrak g_{\one} = \sum_{\lambda}\oplus\, \mathfrak h_{\one}^{\lambda},\label{deoad}\\
\mathfrak g_{\one} = \sum_{\lambda}\oplus\, \mathfrak g_{\one}^{\lambda}\label{deole} 
\end{eqnarray}
Clearly,
\begin{eqnarray}
\mathfrak h_{\one}^{\lambda} = \mathfrak h^{\lambda}, \quad
\mathfrak g_{\one}^{\lambda} = \mathfrak g^{\lambda},\quad \lambda
\neq 0;\qquad \mathfrak h_{\one}^0 = \mathfrak h_{\one}.
\label{sootv} 
\end{eqnarray}
In general,  $\one \not\in \mathfrak g_{\one}^0$.
\begin{lemma}\label{canon}
The decomposition (\ref{decle}) is canonical if and only if $\one
\in \mathfrak g_{\one}^0$.
\end{lemma}
\begin{proof}
If (\ref{decle}) is canonical then $\mathfrak h^0 = \mathfrak
g^0$. Since $\mathfrak h\cdot\one = \mathfrak h = \mathfrak h^0 =
\mathfrak g^0$, we get $\one \in \mathfrak g_{\one}^0$ .

It is sufficient to prove the converse for $\mathbb F=\mathbb C$
(the case $\mathbb F=\mathbb R$ can be easily reduced to it).
Since $L$ is a representation, for all linear functionals
$\lambda,\mu$
\begin{eqnarray}
 \mathfrak h_{\one}^{\lambda}\cdot \mathfrak g_{\one}^{\mu} \subseteq
\mathfrak g_{\one}^{\lambda+\mu}.\label{hinge} 
\end{eqnarray}
If $\one \in \mathfrak g_{\one}^0$ then $\mathfrak
h_{\one}^{\lambda}\subseteq\mathfrak h_{\one}^{\lambda} \cdot
\mathfrak g_{\one}^0$. By (\ref{hinge}), $\mathfrak
h_{\one}^{\lambda} \subseteq \mathfrak g_{\one}^{\lambda}$.
Therefore, each summand in (\ref{deoad}) is contained in the
corresponding summand of (\ref{deole}); since the sums are equal,
the summands coincide. For (\ref{decad}) and (\ref{decle}), the
assertion holds  due to (\ref{sootv}).
\end{proof}
\begin{lemma}\label{soprg}
Let $g \in G_\one$. Then
$$
\mathfrak g = \sum_{\lambda}\oplus\, g(\mathfrak
g_{\one}^{\lambda})
$$
is the decomposition (\ref{deole}) for the Cartan subalgebra
$\Ad(g) \mathfrak h$.
\end{lemma}
\begin{proof}
This is true since $ L(\Ad(g)h) = g L(h) g^{-1}$.
\end{proof}
\begin{lemma}\label{cande}
If the decomposition (\ref{decle}) is canonical then
\begin{eqnarray}
\mathfrak g^{\lambda} \cdot \mathfrak g^{\mu} \subseteq \mathfrak g^{\lambda+\mu}\label{decan} 
\end{eqnarray}
for all $\lambda,\mu$. In particular, $\mathfrak g^0$ is a
subalgebra of the left-symmetric algebra $\mathfrak g$.
\end{lemma}
\begin{proof} The inclusion holds due to (\ref{sootv}),
(\ref{hinge}) and the definition of the canonical decomposition.
\end{proof}
\begin{theorem}\label{maint}
Each complete left-symmetric algebra admits the unique canonical
decomposition.
\end{theorem}
\proof
Let $\mathfrak h$ be any Cartan subalgebra. By (\ref{sootv}),
$\mathfrak g_{\one}^{\lambda} \subseteq \mathfrak g$ if
$\lambda\neq 0$  in (\ref{deole}). Hence
$$
\mathfrak g_{\one}^0 \cap(\one + \mathfrak g) \ne \emptyset.
$$
If $x \in \mathfrak g_{\one}^0 \cap (\one + \mathfrak g) $ then
$g(x) = \one$ for some $g \in G$ (by the definition of the
completeness; recall that the affine action of $G$ may be realized
in the hyperplane $\one+\mathfrak g$). It follows from
Lemma~\ref{soprg} that for the decomposition (\ref{decle}) with
respect to the Cartan subalgebra $\Ad(g) \mathfrak h$ the
condition $\one\in \mathfrak g_{\one}^0$ holds. This makes
possible to apply Lemma~\ref{canon}.

To prove the uniqueness, we use the following known fact: any two
Cartan subalgebras in a (real) solvable Lie algebra are conjugated
by an inner automorphism (if the left symmetric algebra $\frak g$
is complete then the Lie algebra $\frak g$ is solvable; this is
proved, for example, in [2]; see also the introduction). Suppose
that $\mathfrak h$ and $\tilde{\mathfrak h}$ are Cartan
subalgebras corresponding to canonical decompositions
(\ref{decle}). We claim that $\mathfrak h = \tilde{\mathfrak h}$.
There exists $g\in G$ such that $\tilde{\mathfrak h} = \Ad(g)
\mathfrak h$. Then $\tilde{\mathfrak g}^0 = g(\mathfrak g^0)$ by
Lemma~\ref{soprg} and $\one \in \tilde{\mathfrak g}^0 \cap
\mathfrak g^0$ by Lemma~\ref{canon}. Hence $\one,g(\one)\in
\tilde{\mathfrak g}^0$. Since the action of $G$ is free,  the
condition $\one\to g(\one)$ uniquely determines $g\in G$. By
Lemma~\ref{cande} and Corollary~\ref{compl}, $\tilde{\mathfrak
g}^0$ is a complete left symmetric algebra. Hence the subgroup
$\tilde G^0$, corresponding to $\tilde{\mathfrak g}^0$ contains
the unique transformation that sends $\one$ to $g(\one)$.
Therefore, $g\in\tilde G^0$. Thus
$$
\mathfrak h = \Ad (g^{-1}) \tilde{\mathfrak h}=\Ad (g^{-1})
\tilde{\mathfrak g}^0=\tilde{\mathfrak g}^0 =\tilde{\mathfrak h}.
\eqno\square$$
\begin{remark}\label{autom}\rm
Due to the uniqueness, each automorphism of a complete left
symmetric algebra keeps the canonical decomposition. On the other
hand, the existence of the canonical decomposition for a complete
left symmetric algebra which is not nilpotent as a Lie algebra
implies the existence of a non-discrete group of automorphisms.
Indeed, the semisimple parts of the representations $L$ and $\ad$
(they coincide) are differentiations of $\mathfrak g$ according to
(\ref{decan}).
\end{remark}

\section{Graphs of one dimensional canonical decompositions}
In this section, we assume that $\mathbb F=\mathbb C$ and
\begin{eqnarray}
 \dim\mathfrak g^{\lambda}\leq1\label{cadim}
\end{eqnarray}
for all spaces $\mathfrak g^{\lambda}$ in the decomposition
(\ref{decle}). Let $ e_0 \in \mathfrak g^0$, $e_0 \ne 0$. The
spectrum of the linear transformation $L( e_0)$ may be identified
with the set of all nontrivial roots; let us denote it by
$\Lambda$. For $\lambda\in \Lambda$ let $ e_{\lambda}$ be the
eigenvector (for $\lambda = 0$, the chosen vector $ e_0$ satisfies
this ). In the base $\mathcal B = \{ e_{\lambda}\colon\;
\lambda\in\Lambda \}$ the multiplication is defined by relations
\begin{eqnarray}\label{cmula}
e_{\lambda} e_{\mu}=c_{\lambda,\mu} e_{\lambda+\mu}.
\end{eqnarray}
As usual, we assume that $c_{\lambda,\mu}=c_{\mu,\lambda} = 0$ and
$e_{\lambda} = 0$ if $\lambda\not\in \Lambda$.

Each one dimensional canonical decomposition corresponds to graphs
$\Gamma_l$ and $\Gamma_r$ relating to left and right
multiplications. The set of vertices for both graphs is $\Lambda$;
it is natural to assume that 0 is a base point. The graph
$\Gamma_l$ contains the edge from $\lambda$ to $\mu$ if and only
if 
\begin{eqnarray}\label{cunen}
c_{\mu-\lambda,\lambda}\ne 0. 
\end{eqnarray}
For $\Gamma_r$, (\ref{cunen}) is replaced by
\begin{eqnarray}\label{cuner}
 c_{\lambda,\mu-\lambda }\ne 0. 
\end{eqnarray}
By definition, $e_0e_{\lambda} = \lambda e_{\lambda}$. Hence
$c_{0,\lambda} = \lambda$. In particular,
\begin{eqnarray}\label{cuzen}
c_{0,\lambda}\ne 0\quad \mbox{if}\quad \lambda\ne 0. 
\end{eqnarray}
Therefore,
\begin{itemize}
\item[{\tt(l1)}] there is the loop at each nonzero vertex in
$\Gamma_l$; \item[{\tt(r1)}] the vertex 0 in $\Gamma_r$ is joined
with all other vertices.
\end{itemize}
\noindent
Graphs $\Gamma_l$ and $\Gamma_r$ mutually determine each other by
the following procedure. Let $\Gamma_l$ ($\Gamma_r$) contains the
edge $\overline{\lambda\mu}$. Let us remove $\lambda$ to 0 and
join the endpoint of the resulting vector with the endpoint of the
original one.
Doing this with all edges in $\Gamma_l$ ($\Gamma_r$) we get
$\Gamma_r$ ($\Gamma_l$).

It follows from (\ref{cunen}), (\ref{cuner}) that
\begin{itemize}
\item[\tt(l2)] if $\Gamma_l$  contains the edge
$\overline{\lambda\mu}$ then it contains the vertex $\mu-\lambda$,
\end{itemize}
and the same property {\tt(r2)} of $\Gamma_r$. Let us abbreviate
the notation
$$L_{\lambda}= L( e_{\lambda}),\quad R_{\lambda}= R(
e_{\lambda}).$$ All operators $R_{\lambda}$ are nilpotent by
Proposition~\ref{realc}. Hence the assumption (\ref{cadim})
implies
\begin{eqnarray}\label{rzeze}
R_0=0. 
\end{eqnarray}
Therefore,
\begin{itemize}
\item[{\tt(l3)}] there are no edge that comes out of the vertex 0
in $\Gamma_l$; \item[{\tt(r3)}] $\Gamma_r$ contains no loops.
\end{itemize}
As a consequence of (\ref{cuzen}) and (\ref{rzeze}) we get
$$
[\mathfrak g,\mathfrak g] = \sum_{\lambda\ne 0}\oplus\, \mathfrak g^{\lambda}. 
$$
Since the Lie algebra $\mathfrak g$ is solvable, there exists a
base in the linear space $\mathfrak g$ such that $L([\mathfrak
g,\mathfrak g])$ is strictly triangular in it; in particular, this
is true for all $L_{\lambda}$, where $\lambda\ne 0$. Hence the
product $ L_{\lambda_1} L_{\lambda_2} \ldots L_{\lambda_m}$ is
nilpotent if $\lambda_k\neq0$ for some
$k\in\left\{1,\dots,m\right\}$. On the other hand, if
$\lambda_1+\ldots+\lambda_m = 0$ then $ L_{\lambda_1}
L_{\lambda_2} \ldots L_{\lambda_m}$ is diagonal in the base
$\mathcal B$. Therefore, if there is a nonzero summand in this sum
then
\begin{eqnarray}\label{trise}
L_{\lambda_1} \ldots L_{\lambda_m} = 0
\end{eqnarray}
In particular,  $L_{\lambda}L_{-\lambda} = L_{-\lambda}
L_{\lambda}= 0$ for all $\lambda\ne 0$; hence
$$
[L_{\lambda},L_{-\lambda}] = L([e_{\lambda}, e_{-\lambda}]) = 0.
$$
Since $[ e_{\lambda}, e_{-\lambda}] = (c_{\lambda,-\lambda}-
c_{-\lambda,\lambda}) e_0,$, by (\ref{cuzen}),
\begin{eqnarray}\label{elami}
[e_{\lambda}, e_{-\lambda}]  = 0.
\end{eqnarray}
Thus $c_{\lambda,-\lambda} = c_{-\lambda,\lambda}$. This means
that
\begin{itemize}
\item[\tt(l4)] if some edge  in $\Gamma_l$ comes to 0 from
$\lambda$ then the same is true for the vertex $-\lambda$ (in
particular, $\Gamma_l$ contains this vertex).
\end{itemize}
Clearly, the analogous property {\tt(r4)} holds for $\Gamma_r$.
Furthermore, (\ref{trise}) immediately implies that
\begin{itemize} \item[{\tt (l5)}] $\Gamma_l$ has no cycle.
\end{itemize}
Suppose that $\Gamma_l$ contains two consecutive edges, say,
$\overline{\lambda\mu}$ and $\overline{\mu\nu}$. Then
$$
e_{\nu-\mu} ( e_{\mu-\lambda} e_{\lambda})\ne 0.
$$
By (\ref{curvi}),
$$
  e_{\nu-\mu} (  e_{\mu-\lambda}   e_{\lambda}) =
  e_{\mu-\lambda} (  e_{\nu-\mu}   e_{\lambda}) +
[  e_{\nu-\mu},   e_{\mu-\lambda}]   e_{\lambda}.
$$
Since at least one summand in the right side above is nontrivial,
we get
\begin{itemize}\item[{\tt(l6)}] each pair of consecutive edges
in $\Gamma_l$ can be included either to a triangle or to a
parallelogram  in $\Gamma_l$.
\end{itemize}
The parallelogram may be degenerate but new edges coincides with
initial ones only if $\mu-\lambda=\nu-\mu$.

It follows from  (\ref{lrcpm}) and Proposition~\ref{realc} that
\begin{eqnarray}\label{trrrr}
 \Tr(R(x)R(y) ) = 0
\end{eqnarray}
for all $x,y\in\mathfrak g$. Analogously, the equality
$$[L(y),R(x)^2]=[L(y),R(x)]R(x)+R(x)[L(y),R(x)],$$
taken together with (\ref{lrcpm}) and (\ref{trrrr}) implies that
$$
\Tr (R(x)^2 R(y))=0.
$$
In particular, $\Tr (R_{\lambda}^2 R_{-2\lambda}) = 0$ for all
$\lambda\in\Lambda$. Therefore,
\begin{itemize}
\item[{\tt (r5)}] $\Gamma_r$ cannot contain exactly one path of
the type
\begin{eqnarray*}\label{putod}
 \mu\to\mu-2\lambda \to \mu-\lambda \to\mu.
\end{eqnarray*}
\end{itemize}
We formulate also several properties that hold only for simple
left symmetric algebras. A subspace $\frak j\subseteq\frak g$ is
an {\it ideal} in the left symmetric algebra $\frak g$ if
$\frak{jg}\subseteq\frak j$ and $\frak{gj}\subseteq\frak j$
(i.\,e. if $\frak j$ is a two-side ideal); an algebra is {\it
simple} if it contains no proper ideals.
\begin{itemize}
\item[{\tt (s1)}] For any $\lambda\in\Lambda$,  the graph
$\Gamma_l$ contains an edge parallel to $\overline{0\lambda}$;
the graph $\Gamma_r$ has an edge that comes out of $\lambda$.
\end{itemize}
These properties are mutually dual with respect to the procedure
described above. They hold since the kernel of the representation
$L$ is an ideal in $\mathfrak g$.
\begin{itemize}
\item[\tt (s2)] The vertex 0 is reachable from any point
$\lambda\in\Lambda$ in the union of graphs $\Gamma_l$ and
$\Gamma_r$.
\end{itemize}
The set $\cal R$ of all vertices that are reachable from $\lambda$
defines a nontrivial  ideal which must coincide with $\mathfrak
g$; hence $\cal R$  contains $e_0$.
\begin{itemize}
\item[{\tt (s3)}] $\Gamma_l$  and $\Gamma_r$ contain at least one
pair of symmetric vertices $\lambda,-\lambda$ and edges from them
to 0.
\end{itemize}
The assertion follows from {\tt(s2), (l4)} and {\tt(r4)}.

The properties above distinguish  graphs of complex simple
algebras for dimensions ${}\le 5$. Clearly, graphs which can be
identified by rotations and dilations are equivalent (these
operations corresponds to the multiplication of $e_0$ by a complex
number). Thus we may assume that $1$ is a vertex and that
$\Gamma_l$ contains no vertex $\lambda\neq0$ such that
$|\lambda|<1$. Also, we drop trivial edges defined in {\tt(l1)}.

There is no simple algebras of dimension 2.

In dimension 3, there is exactly one graph that is determined by
{\tt(s3)}. The corresponding algebra is unique up to an
isomorphism and can be defined by relations
$$
  e_0   e_1 =   e_1,\quad
  e_0   e_{-1} = -  e_{-1},\quad
  e_1   e_{-1} =   e_{-1}   e_1 =   e_0
$$
(Auslander's algebra; it defines an affine action that gives an
example of a free transitive group of affine motions that contains
no translation, see [3]).

In dimension 4, there is also only one graph $\Gamma_l$
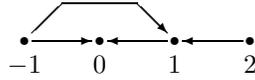
\begin{center}
\unitlength=1mm \linethickness{0.4pt}
\begin{picture}(31.00,10.00)
\put(1.00,5.00){\vector(1,0){8.00}}
\put(19.00,5.00){\vector(-1,0){8.00}}
\put(29.00,5.00){\vector(-1,0){8.00}}
\put(1.00,6.00){\line(1,1){4.00}}
\put(5.00,10.00){\line(1,0){10.00}}
\put(15.00,10.00){\vector(1,-1){4.00}}
\put(0.00,3.00){\makebox(0,0)[ct]{$-1$}}
\put(10.00,3.00){\makebox(0,0)[ct]{$0$}}
\put(20.00,3.00){\makebox(0,0)[ct]{$1$}}
\put(30.00,3.00){\makebox(0,0)[ct]{$2$}}
\put(0.00,5.00){\circle*{1.00}} \put(10.00,5.00){\circle*{1.00}}
\put(20.00,5.00){\circle*{1.00}} \put(30.00,5.00){\circle*{1.00}}
\end{picture}
\end{center}
and the unique up to an isomorphism simple algebra
\begin{eqnarray*}
&  e_0   e_{-1} = -  e_{-1},\quad
  e_0   e_1 =   e_1,\quad
  e_0   e_2 = 2  e_2,\\
&  e_{-1}   e_1 =   e_1   e_{-1} =   e_0,\quad
  e_2   e_{-1} = 2  e_1,\quad
  e_{-1}   e_2 =   e_1 .
\end{eqnarray*}
In dimension 5, there is one parameter family of graphs and simple
algebras. Graphs consists of two pairs of symmetric vertices and
edges coming to 0. Algebras can be defined by relations
\begin{eqnarray*}
&  e_0   e_{-1} = -  e_{-1},\quad
  e_0   e_1 =   e_1,\quad
  e_{-1}   e_1 =   e_1   e_{-1} =   e_0,\\
&  e_0   e_{-\lambda} = -\lambda   e_{-\lambda},\quad
  e_0   e_{\lambda} = \lambda   e_{\lambda},\quad
  e_{-\lambda }   e_{\lambda} =   e_{\lambda}   e_{-\lambda } =   e_0.
\end{eqnarray*}
For $\lambda = 2$, some additional relations may hold:
\bigskip\par
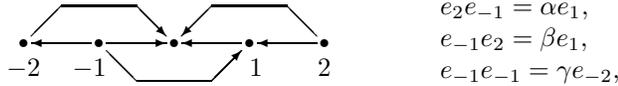
\parbox{0.4\textwidth}{\unitlength=1mm
\linethickness{0.4pt}
\begin{picture}(41.00,10.00)
\put(9.00,5.00){\vector(-1,0){8.00}}
\put(11.00,5.00){\vector(1,0){8.00}}
\put(29.00,5.00){\vector(-1,0){8.00}}
\put(39.00,5.00){\vector(-1,0){8.00}}
\put(1.00,6.00){\line(1,1){4.00}}
\put(5.00,10.00){\line(1,0){10.00}}
\put(15.00,10.00){\vector(1,-1){4.00}}
\put(39.00,6.00){\line(-1,1){4.00}}
\put(35.00,10.00){\line(-1,0){10.00}}
\put(25.00,10.00){\vector(-1,-1){4.00}}
\put(11.00,4.00){\line(1,-1){4.00}}
\put(15.00,0.00){\line(1,0){10.00}}
\put(25.00,0.00){\line(1,1){4.00}}
\put(25.00,0.00){\vector(1,1){4.00}}
\put(0.00,5.00){\circle*{1.00}} \put(10.00,5.00){\circle*{1.00}}
\put(20.00,5.00){\circle*{1.00}} \put(30.00,5.00){\circle*{1.00}}
\put(40.00,5.00){\circle*{1.00}}
\put(0.00,3.00){\makebox(0,0)[ct]{$-2$}}
\put(11.00,3.00){\makebox(0,0)[rt]{$-1$}}
\put(30.00,3.00){\makebox(0,0)[lt]{$1$}}
\put(40.00,3.00){\makebox(0,0)[ct]{$2$}}
\end{picture}
}\hfill
\parbox{0.5\textwidth}{ $   e_2   e_{-1} =
\alpha   e_1$,\\ $  e_{-1}   e_2 = \beta   e_1$,\\
$  e_{-1}  e_{-1} = \gamma   e_{-2}$,}

\bigskip

\noindent where $2\alpha = \beta+\gamma$. For $\alpha=0$,
$\beta=0$, $\gamma=0$ edges $-1\to 1$, $2\to 1$, $-1\to -2$ must
be dropped, respectively. Algebras corresponding to triples
$(\alpha,\beta,\gamma)$ и $(\alpha',\beta',\gamma')$ are
isomorphic if and only if vectors $(\alpha,\beta,\gamma)$ and
$(\alpha',\beta',\gamma')$ are collinear (see Remark~\ref{autom}).
Thus the family corresponding to $\lambda = 2$ is modelled by the
projective line $\mathbb F \mathbb P^1$ with three distinguished
points.

There is only one simple algebra with one dimensional root
decomposition which is not mentioned above:
\bigskip\par

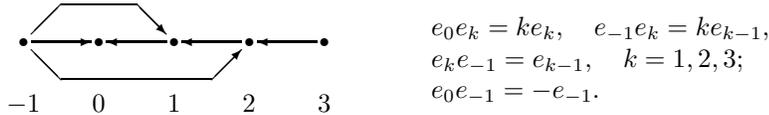
\parbox{0.4\textwidth}{\unitlength=1mm
\linethickness{0.4pt}
\begin{picture}(41.00,15.00)
\put(1.00,10.00){\vector(1,0){8.00}}
\put(19.00,10.00){\vector(-1,0){8.00}}
\put(29.00,10.00){\vector(-1,0){8.00}}
\put(39.00,10.00){\vector(-1,0){8.00}}
\put(1.00,11.00){\line(1,1){4.00}}
\put(5.00,15.00){\line(1,0){10.00}}
\put(15.00,15.00){\vector(1,-1){4.00}}
\put(1.00,9.00){\line(1,-1){4.00}}
\put(5.00,5.00){\line(1,0){20.00}}
\put(25.00,5.00){\vector(1,1){4.00}}
\put(0.00,10.00){\circle*{1.00}} \put(10.00,10.00){\circle*{1.00}}
\put(20.00,10.00){\circle*{1.00}} \put(30.00,10.00){\circle*{1.00}}
\put(40.00,10.00){\circle*{1.00}}
\put(0.00,3.00){\makebox(0,0)[ct]{$-1$}}
\put(10.00,3.00){\makebox(0,0)[ct]{$0$}}
\put(20.00,3.00){\makebox(0,0)[ct]{$1$}}
\put(30.00,3.00){\makebox(0,0)[ct]{$2$}}
\put(40.00,3.00){\makebox(0,0)[ct]{$3$}}
\end{picture}
}\hfill
\parbox{0.5\textwidth}{$  e_0   e_k = k  e_k,\quad
  e_{-1}   e_k = k  e_{k-1}$,\\ $  e_k
  e_{-1} =   e_{k-1},\quad k=1,2,3$;\\ $  e_0
  e_{-1} = -  e_{-1} $.}
\bigskip

\noindent It is included to infinite series and an infinite
dimensional left symmetric algebra with the same relations and $k$
running over $\mathbb N$.

\section{Acknowledgments}
I am  grateful to I.G.\,Korepanov who encouraged me to prepare the
paper for arXiv and to D.\,Burde who translated it to German and
used in his paper \cite{Bu}.

\end{document}